\documentclass[12pt]{article}
\usepackage{amsmath,amsopn,amssymb,amsthm,amsfonts}
\usepackage[T2A]{fontenc}
\usepackage[cp1251]{inputenc}
\usepackage{longtable}

\newtheorem{theorem}{Theorem}[section]

\newtheorem{definition}[theorem]{Definition}

\newcommand{\w}{\omega}
\newcommand\g{{\mathfrak g}}
\newcommand\G{{\mathfrak G}}
\newcommand\h{{\mathfrak h}}

\newcommand\ag{{\mathfrak a}}

\newcommand\n{{\mathfrak n}}

\begin{document}

{\bf \large
\centerline{Smolentsev~N.~K., Chernova~K.~V.}

\vspace{3mm}
\centerline{On left-invariant semi-K\"{a}hler structures}
\centerline{on six-dimensional nilpotent nonsymplectic Lie groups}}
\vspace{3mm}

\begin{abstract}
It is known that there are 34 classes of six-dimensional nilpotent Lie groups, many of which admit left-invariant symplectic and complex structures. Among them there are three classes of groups on which there are no left-invariant symplectic structures, but there are complex structures. In this work, semi-K\"{a}hler and almost para-semi-K\"{a}hler structures are defined on such Lie groups and their geometric properties are studied.
\end{abstract}

\section{Introduction} \label{Intro}
It was shown in \cite{Goze-Khakim} that 26 out of 34 classes of six-dimensional nilpotent Lie groups admit left-invariant symplectic structures. Of the remaining eight classes of nonsymplectic Lie groups, five Lie groups also do not admit complex structures and three nonsymplectic Lie groups admit complex structures \cite{Sal}. The first 5 groups were studied in the work of Smolentsev N.K. \cite{Smolen-19}. The geometry of the last three nonsymplectic Lie groups is currently insufficiently studied.
The goal of the work is to determine new left-invariant geometric structures on the indicated three six-dimensional Lie groups, compensating in a sense for the absence of symplectic structures, as well as to study their geometric properties.

Instead of the symplecticity property of the 2-form $\w$, in this work we take the semi-K\"{a}hler property $\omega\wedge d\omega = 0$ on a six-dimensional manifold.

In this paper, families of semi-K\"{a}hler 2-forms $\omega$, compatible complex structures and corresponding pseudo-Riemannian metrics are found for each group.
The properties of left-invariant 2-forms $\omega$ having a non-degenerate exterior derivative $d\omega$ in the sense of Hitchin are determined and studied. The geometric properties of para-semi-K\"{a}hler and semi-K\"{a}hler structures on the indicated Lie groups are studied. The associated metrics are determined and their curvature properties are calculated.

\section{Preliminaries} \label{Pre}
Let us recall the basic concepts used in this work.

Let $G$ be a real Lie group and $\g$ its Lie algebra.
A Lie group is called nilpotent if its Lie algebra is nilpotent, that is, the following decreasing sequence of ideals $\g_k$ is finite: $\g_0 = \g, \ \g_{k+1} = [\g, \g_k]$.

An almost complex structure on a $2n$-dimensional manifold $M$ is a field $J$ of endomorphisms $J: TM\to TM$ such that $J^2 =-Id$.
The almost complex structure $J$ is called integrable if the Nijenhuis tensor, defined by the equality $N_J(X,Y) = [JX,JY ] -[X,Y] -J[JX,Y] -J[X,JY]$, vanishes, for any vector fields $X, Y$ on $M$. In this case, $J$ defines the structure of a complex manifold on $M$.

We will consider left-invariant almost complex structures on the Lie group $G$, which are definite by the left-invariant endomorphism field $J: TG \to TG$ of the tangent bundle $TG$.
Since such a tensor $J$ is determined by a linear operator on the Lie algebra $\g = T_eG$, we will say that $J$ is an invariant almost complex structure on the Lie algebra $\g$.
In this case, the integrability condition for $J$ is formulated at the Lie algebra level: $N_J(X,Y) = 0$, for any $X,Y\in \g$.
In this case, we will say that $J$ is a complex structure on the Lie algebra $\g$.

A left-invariant symplectic structure $\omega$ on a Lie group $G$ is definite by a closed 2-form of maximum rank. The closedness of the form $\omega$ is equivalent to the following condition on the Lie algebra $\g$: $\omega([X,Y],Z) -\omega([X,Z],Y) + \omega([Y,Z],X) = 0$, $\forall X,Y,Z\in \g$. In this case, the Lie algebra $\g$ will be called symplectic.

A left-invariant K\"{a}hler structure on a Lie group $G$ is a triple $(\w, J, g)$ consisting of a left-invariant symplectic form $\w$, a left-invariant complex structure $J$ compatible with the form $\w$: $\w(JX,JY) = \w(X,Y)$, and the (pseudo)Riemannian metric $g$, defined by the formula $g(X,Y) = \w(X,JY)$.

If the symplectic form $\w$ is given, then to find a K\"{a}hler structure on a Lie group it is necessary to find a complex structure $J$ compatible with $\w$.
To do this, it is sufficient to solve the following system of equations: $\w(JX,Y) + \w(X,JY) = 0$, $J^2 = -Id$ and $N_J = 0$.
If the symplectic form $\w$ and operator $J$ have matrices $\w_{ij}$ and $J_j^i$ in the basis $\{e_1,\dots , e_{2n}\}$ of the Lie algebra.
Then the system of equations for finding the K\"{a}hler structure $(\w, J)$ consists of the following equations for the variables $J_j^i$:
\begin{equation} \label{0}
\left\{
\begin{tabular}{l}
$\w_{kj}J^k_i+\w_{ik}J^k_j=0$, \\
$J_k^i\, J^k_j = -\delta_j^i,$  \\
$J_i^l J_j^m C_{lm}^k-J_i^l J_m^k C_{lj}^m-J_j^l J_m^k C_{il}^m -C_{ij}^k=0$,
 \end{tabular}
\right.
\end{equation}
where $\delta_j^i$ is the identity matrix, $C_{ij}^k$ are the structure constants of the Lie algebra, and the indices vary from 1 to $2n$.

Let an almost complex structure $J$ be given on a Lie algebra $\g$. Let us define an increasing sequence of $J$-invariant ideals:
$$
\mathfrak{a}_0(J) = 0, \ \mathfrak{a}_s(J) = \{X \in \g | [X, \g] \subset \mathfrak{a}_{s-1}(J) \mbox{ and } [JX, \g] \subset \mathfrak{a}_{s-1} (J)\}, \ s \ge 1.
$$
In particular, the ideal $\mathfrak{a}_1(J)$ lies in the center and has dimension at least two.

\begin{definition}\label{1}
A leftinvariant almost complex structure $J$ on a Lie group $G$ is called nilpotent if for a series of ideals $\mathfrak{a}_s(J)$ defined above there is a number $p$ such that $\mathfrak{a}_p(J) = \g$.
\end{definition}

It was shown in \cite{CFU} that K\"{a}hler structures on six-dimensional nilpotent Lie groups have nilpotent complex structures.

An almost para-complex structure on a $2n$-dimensional manifold $M$ is the field $P$ of endomorphisms of the tangent bundle $TM$ such that $P^2 = Id$, and the ranks of the eigendistributions $T^\pm M : = ker(Id\mp P)$ are equal.
An almost paracomplex structure $P$ is said to be integrable if the distributions $T^\pm M$ are involutive. In this case, $P$ is called a paracomplex structure.

In the case of a left-invariant paracomplex structure $P$ on a Lie group $G$, the invariant subspaces $T_e^\pm G$ are subalgebras.
Therefore, the paracomplex Lie algebra $\g$ can be represented as a direct sum of two subalgebras:
$$
\g = \g^+ \oplus \g^-.
$$
The left-invariant almost paracomplex structure $P$ defines a paracomplex structure on the Lie group $G$ in the case when the Nijenhuis tensor vanishes:
$$
N_P(X,Y) = [PX,PY] + [X,Y] - P[PX,Y] - P[X,PY] = 0, \forall X,Y \in \g.
$$
A review of the theory of para-complex structures is presented in \cite{Aleks}.

A left-invariant para-K\"{a}hler structure on a Lie group $G$ is a triple $(\w,P,g)$ consisting of a left-invariant symplectic form $\w$, a left-invariant paracomplex structure $P$ compatible with the form $\w$: $\w(PX,PY) =-\w(X,Y)$, and the (pseudo)Riemannian metric $g$, defined by the formula $g(X,Y) = \w(X,PY)$.

Note that the term “para-K\"{a}hler” manifold is also used in a different sense. Namely, almost Hermitian manifolds satisfying the K\"{a}hler identity $g(R(X,Y)Z,W) = g(R(X,Y)JZ,JW)$, where $R$ is the curvature tensor and $J$ is an almost complex structure compatible with Riemannian metric $g$ are also called para-K\"{a}hler manifolds.
In this paper, para-K\"{a}hler manifolds are considered from the point of view of paracomplex geometry, which was introduced by P.~Liberman in 1952 by analogy with complex geometry.
Note for comparison that such para-K\"{a}hler manifolds satisfy the following identity: $g(R(X,Y)Z,W) = -g(R(X,Y)PZ, PW)$, where $P$ is the tensor of the almost paracomplex structure compatible with the pseudo -Riemannian metric $g$.

For the nonsymplectic Lie groups under consideration, the condition $d\w = 0$ leads to degeneracy of the 2-form $\w$.
At the same time, $d(\w^3) = 3\w\wedge\w\wedge d\w = 0$ for any 2-form $\w$ on a six-dimensional manifold.
An intermediate property would be $d(\w^2) = 2\w\wedge d\w = 0$.
Therefore, we reduce the closedness property $d\w = 0$ and require that the following property hold:
$$
\w\wedge d\w = 0.
$$
Note that in the case of a left-invariant almost Hermitian structure on a Lie group $G$ of dimension $2n$, the property $d(\w^{n-1}) = 0$ of the fundamental form $\w$ defines a class of semi-K\"{a}hler manifolds according to the Gray–Harvella classification \cite{Gray-Harv}, i.e. such that $\delta \w = 0$.
In our case of pseudo-Hermitian metrics on Lie groups, we will by analogy call such manifolds semi-K\"{a}hler.

\begin{definition}\label{1}
A left-invariant almost Hermitian (almost para-Hermitian) structure $(G,g,\w,J)$ of dimension $2n$, whose fundamental form $\w$ has the property $d(\w^{n-1}) = 0$ is called semi-K\"{a}hler (semi-para-K\"{a}hler).
\end{definition}

If the non-degenerate 2-form $\w$ is not closed, then we can consider the 3-form $d\w$.
In \cite{Hitch}, Hitchin defined the concept of non-degeneracy (stability) for 3-forms $\Omega$ and constructed a linear operator $K_\Omega$ whose square is proportional to the identity operator $Id$. Let us recall Hitchin's construction.
Let $V$ be a 6-dimensional real vector space and $\mu$ a volume form on $V$. Let $\Omega\in \Lambda^3 V^*$ and $X\in V$, then $i_X\Omega\in \Lambda^2 V^*$ and $i_X\Omega\wedge\Omega\in  \Lambda^5 V^*$.
The natural pairing by the exterior product $V^*\otimes \Lambda^5 V^* \to \Lambda^6 V^* \cong \mathbb{R}\mu$ defines an isomorphism $A:\Lambda^5 V^* \cong  V$ and using this, Hitchin defined the linear transformation $K_\Omega : V \to V$ by the following formula:
$$
K_\Omega(X) = A(i_X\Omega\wedge\Omega).
$$
In other words, $i_{K_\Omega(X)}\mu = i_X\Omega\wedge\Omega$.
The operator $K_\Omega$ has the following properties: $\mathrm{trace}(K_\Omega) = 0$ and $K_\Omega^2 = \lambda(\Omega)Id$. Therefore, in the case $\lambda(\Omega)< 0$, we obtain the structure $J_\Omega$ of a complex vector space on the space $V$:
$$
J_\Omega=\frac{1}{\sqrt{-\lambda(\Omega)}} K_\Omega,
$$
and if $\lambda(\Omega)> 0$, then we obtain the para-complex structure $P_\Omega$, i.e., $P_\Omega^2 = Id$, $P_\Omega\ne 1$ by a similar formula:
$$
P_\Omega=\frac{1}{\sqrt{\lambda(\Omega)}} K_\Omega.
$$
Thus, the operator $K_{d\w}$ can define either an almost complex or almost paracomplex structure on a 6-dimensional Lie group when $d\w$ is non-degenerate and $\mu = \w^n$.

Let $\nabla$ be the Levi-Civita connection corresponding to the (pseudo)Riemannian metric $g$.
It is determined from a six-term formula \cite{KN}, which for left-invariant vector fields $X,Y,Z$ on a Lie group takes the form: $2g(\nabla_XY, Z) = g([X,Y],Z) + g([Z,X],Y ) + g(X,[Z,Y])$.
The curvature tensor is defined by the formula $R(X,Y) = [\nabla_X, \nabla_Y] - \nabla_{[X,Y]}$.
The Ricci tensor $Ric(X,Y)$ is the convolution of the curvature tensor along the first and fourth (upper) indices.
The Ricci operator $RIC$ is defined by the formula $Ric(X,Y) = g(RIC(X),Y)$.
To calculate the geometric characteristics of left-invariant structures on Lie groups, the Maple system was used.

\section{Nonsymplectic Lie groups} \label{Groups}
In this paper we study nilpotent Lie groups $G_i$, $i=1,2,3$, which do not admit left-invariant symplectic structures, but can have complex structures. In accordance with the classification \cite{Sal}, they have Lie algebras with the following non-zero Lie brackets of basis vectors $e_1,\dots, e_6 \in \g_i$:

$\g_1: [e_1,e_2] =e_4, \ [e_2,e_3] =e_5, \ [e_1,e_4] =e_6, \ [e_3,e_5] =-e_6.$

$\g_2: [e_1,e_2]=e_4, \ [e_1,e_4] =e_5, \ [e_2, e_4] = e_6.$

$\g_3: [e_1, e_2] = e_6, \ [e_3, e_4] = e_6.$

\subsection{Lie group $G_1$} \label{G1}
Let us consider the first Lie group, which does not admit left-invariant symplectic structures, but has complex ones. Complex structures on this Lie group were studied in the work of Magnin \cite{Magn}.
We will use the following commutation relations of the Lie algebra $\g_1$, corresponding to \cite{Magn} (which are obtained by replacing: $e_2 \leftrightarrow e_3, \ e_5 = -e_5$):
$$
[e_1,e_3] =e_4, \ [e_1,e_4] =e_6, \ [e_2,e_3] =e_5, \ [e_2,e_5] =e_6.
$$

\subsubsection{Semi-para-K\"{a}hler structures.}
Let $\w = w_{ij}e^i\wedge e^j$ be an arbitrary 2-form in the dual basis $\{e^i\}$.
%Let us consider the question of the para-K\"{a}hler property and non-degeneracy in the Hitchin sense of its exterior derivative $d\w$ and study the paracomplex structures $P_{d\w}$ corresponding to non-degenerate $d\w$.
The property $\w\wedge d\w = 0$ is satisfied under the following conditions:

\centerline{$-w_{24}w_{56} + w_{25}w_{46} - w_{26}w_{45} = 0, \ -w_{14}w_{56} + w_{15}w_{46} - w_{16}w_{45} = 0,$}

\centerline{$w_{13}w_{46}- w_{14}w_{36} + w_{16}w_{34} + w_{23}w_{56} - w_{25}w_{36} + w_{26}w_{35} = 0.$}

There are 2 solutions $\w_1$ and $\w_2$ of this system of equations when the 3-form $d\w$ is non-degenerate:
$$
\w_1=e^1\wedge\left(w_{12} e^2+w_{13} e^3+ \frac{w_{15} w_{46} -w_{16} w_{45}}{w_{56}}\, e^4+w_{15} e^5 +w_{16} e^6 \right)+$$
$$
+e^2\wedge\left(\Omega_{23} e^3+ \frac{w_{25} w_{46} -w_{26} w_{45}}{w_{56}} \, e^4+w_{25} e^5+w_{26} e^6\right)+$$
$$
+e^3\wedge\left(w_{34} e^4 +w_{35} e^5 +w_{36} e^6 \right)+e^4\wedge\left(w_{45} e^5+w_{46} e^6 \right)+w_{56} e^5\wedge e^6,$$
where
$$
\Omega_{23}=-\frac{w_{13} w_{46} w_{56} -w_{15} w_{36} w_{46} +w_{16} w_{34} w_{56} +w_{16} w_{36} w_{45} -w_{25} w_{36} w_{56} +w_{26} w_{35} w_{56}}{w_{56}^2}
$$
and
$$
\w_2=e^1\wedge\left(w_{12} e^2+\Omega_{13} e^3 +w_{14} e^4 +\frac{w_{16} w_{45}}{w_{46}}\, e^5 +w_{16} e^6 \right)+$$
$$
+e^2\wedge\left(w_{23} e^3+ w_{24} e^4+ \frac{w_{26} w_{45}}{w_{46}}\, e^5 +w_{26} e^6 \right)+ $$
$$
+e^3\wedge\left(w_{34} e^4 +w_{35} e^5 +w_{36} e^6 \right) +e^4\wedge(w_{45} e^5 +w_{46} e^6 ),$$
where
$$
\Omega_{13}=\frac{w_{14} w_{36} w_{46} -w_{16} w_{34} w_{46} -w_{26} w_{35} w_{46} +w_{26} w_{36} w_{45}}{w_{46}^2}.
$$

Let us present the second solution in more detail.
Calculations show that for a given 2-form $\w_2$ the exterior derivative $d\w_2$ is non-degenerate in the sense of Hitchin.
In this case, $\lambda(d\w) = w_{46}^4 > 0$.
Therefore, the 3-form $d\w_2$ defines an almost paracomplex structure $P_{d\w_2}$ on the space $\g_1$ as follows:
$$
P_{d\w_2}=\frac{1}{\sqrt{\lambda(d\w_2)}} K_{d\w_2}.
$$
Matrix expression the $P_{d\w_2}$:
\begin{equation} \label{1}
P_{d\w_2} =\left[
  \begin{array}{cccccc}
    1 &0&0&0&0&0 \\
    0& -1 &0&0&0&0 \\
    0&0& 1 &0&0&0 \\
    \frac{-2w_{16}}{w_{46}} &0& \frac{-2w_{36}}{w_{46}} & -1 &0&0 \\
    \frac{-2w_{26}}{w_{46}} &0& \frac{-2w_{45}}{w_{46}} &0& -1&0 \\
    \frac{2w_{16}w_{36}+2w_{26}w_{45}}{w_{46}^2}&\frac{-2w_{24}w_{46} +2w_{26}w_{36}}{w_{46}^2} & \frac{2w_{36}^2+2w_{45}^2}{w_{46}^2} & \frac{2w_{36}}{w_{46}} & \frac{2w_{45}}{w_{46}}&1\\
    \end{array}
\right].
\end{equation}

The almost paracomplex structure $P_{d\w_2}$ is not integrable; the Nijenhuis tensor does not vanish.
For an operator $P_{d\w_2}$ of almost paracomplex structure, the property of compatible with the form $\w_2$ is satisfied:
$$
\w_2(P_{d\w_2}X,P_{d\w_2}Y) = -\w_2(X,Y).
$$

Let's define the associated metric $g_{d\w_2}(X,Y) = \w_2(X,P_{d\w_2}Y)$. We obtain a family of semi-para-K\"{a}hler structures depending on 11 parameters $w_{ij}$.
Calculations show that the associated metric has non-zero scalar curvature:
$$
S(\w_2) =-\frac{w_{46}^2}{(w_{35}w_{46} -w_{36}w_{45})(w_{12}w_{46} -w_{14}w_{26} +w_{16}w_{24})}.
$$

When calculating the geometric characteristics of the pseudo-Riemannian metric $g_{d\w_2}$, very complicated expressions are obtained.
Therefore, we will consider a special case when some of the parameters $w_{ij}$ are equal to zero.
First, note that the tensor $P_{d\w_2}$ and the scalar curvature $R(\w_2)$ do not depend on the parameters $w_{23}$ and $w_{34}$.
We will consider them zero. In addition, from the requirement that the maximum number of components of the Nijenhuis tensor be equal to zero, we obtain: $w_{26} = 0$, $w_{36} = 0$, $w_{45} = 0$.
In this case, we have a simpler expression for $\w_2$ and $P_{d\w_2}$:

$$
\w_2 =\left[
  \begin{array}{cccccc}
    0&w_{12}&0&w_{14}&0&w_{16} \\
    -w_{12}& 0&0&w_{24}&0&0 \\
    0&0& 0 &0&w_{35}&0 \\
    -w_{14}&-w_{24}& 0 & 0&0&w_{46} \\
    0&0& -w_{35}&0& 0&0 \\
    -w_{16}&0&0&-w_{46}& 0&0\\
    \end{array}
\right], \quad
P_{d\w_2} =\left[
  \begin{array}{cccccc}
    1 &0&0&0&0&0 \\
    0& -1 &0&0&0&0 \\
    0&0& 1 &0&0&0 \\
    \frac{-2w_{16}}{w_{46}} &0& 0 & -1 &0&0 \\
    0&0& 0&0& -1&0 \\
    0&\frac{-2w_{24}}{w_{46}}& 0& 0& 0&1\\
    \end{array}
\right].
$$
Then the metric tensor $g_{d\w_2}$ has the following Ricci tensor:

$$
Ric =\frac{w_{46}^2}{2w_{35}(w_{12}w_{46}+w_{16}w_{24})}\left[
  \begin{array}{cccccc}
   \frac{-2w_{14}w_{16}}{w_{46}}&w_{12}&0&-w_{14}&0&w_{16} \\
    w_{12}& 0&0&-w_{24}&0&0 \\
    0&0& 0 &0&w_{35}&0 \\
    -w_{14}&-w_{24}& 0 & 0&0&w_{46} \\
    0&0& w_{35}&0& 0&0 \\
    w_{16}&0&0&w_{46}& 0&0\\
    \end{array}
\right]
$$
and scalar curvature $S=\frac{w_{46}^2}{w_{35}(w_{12}w_{46}+w_{16}w_{24})}.$

\subsubsection{Semi-K\"{a}hler structures.}
In this section, we consider the complex structure on $\g_1$ found in Magnin’s work \cite{Magn}:

\begin{equation} \label{2}
J =\left[
  \begin{array}{cccccc}
    0& -1 &0&0&0&0 \\
    1& 0  &0&0&0&0 \\
    0&0& 0 &0&0&\xi_{36} \\
    0&0&0 & 0& -1 &0 \\
    0&0& 0& 1 & 0&0 \\
    0&0& -\xi_{36} & 0& 0&0\\
    \end{array}
\right],
\end{equation}
where $\xi_{36} =\pm 1$.
The remaining complex structures on $\g_1$ are obtained by the action of the group of automorphisms.
For definiteness, let us accept $\xi_{36} =1$.

Let us consider 2-forms $\w = w_{ij} e^i\wedge e^j$, compatible with the operator $J$: $\w(JX,JY) = \w(X,Y)$, and having the semi-K\"{a}hler property $\w \wedge d\w = 0$. There are 4 types of such 2-forms:

\begin{multline*}
\w_1 =e^1\wedge \left(w_{12} e^2 +w_{13} e^3 -\frac{w_{13}w_{45}}{w_{34}}   e^5 +\frac{w_{13}w_{35}}{w_{34}} e^6 \right) +e^2\wedge \left(\frac{w_{13}w_{35}}{w_{34}}e^3+\frac{w_{13}w_{45}}{w_{34}}e^4 -w_{13}e^6 \right)+
\\
e^3\wedge (w_{34}e^4 +w_{35}e^5 +w_{36}e^6 ) +e^4\wedge (w_{45} e^5 -w_{35} e^6) +w_{34}e^5\wedge e^6,
\end{multline*}
\begin{multline*}
\w_2 =e^1\wedge\left(w_{12}e^2 -\frac{w_{23}w_{45}}{w_{35}}e^5 +w_{23}e^6\right)
+e^2\wedge\left(w_{23}e^3 +\frac{w_{23}w_{45}}{w_{35}}e^4 \right)+ \\
e^3\wedge(w_{35}e^5 +w_{36} e^6)
+e^4\wedge(w_{45}e^5 -w_{35}e^6),
\end{multline*}
\begin{multline*}
\w_3 =e^1\wedge\left(w_{12} e^2 -\frac{w_{15}w_{34}}{w_{45}}e^3 +w_{15}e^5 \right)
+e^2\wedge\left(-w_{15}e^4 +\frac{w_{15}w_{34}}{w_{45}}e^6\right)+\\
e^3\wedge(w_{34}e^4 +w_{36}e^6 )
+w_{45}e^4\wedge e^5 +w_{34}e^5\wedge e^6,
\end{multline*}
\begin{multline*}
\w_4=e^1\wedge(w_{12}e^2 +w_{13}e^3) -w_{13}e^2\wedge e^6 +e^3\wedge(w_{34}e^4 +w_{36}e^6) +w_{34}e^5\wedge e^6.
\end{multline*}

As a result, we obtain four families of semi-K\"{a}hler structures $(J,\w_k,g_{k})$, where $g_{k}(X,Y) =\w_k (X,JY)$, $k=1, 2, 3, 4$.
For each of them, the geometric characteristics of the associated metric $g_{k}$ are easily calculated.

Let us consider the last solution $\w_4$ in more detail.
The associated pseudo-Riemannian metric $g_{4}(X,Y) =\w_4 (X,JY)$ has the following Ricci tensor
$$
Ric_4 =\frac{1}{w_{34}^2 w_{12}}
\left[
  \begin{array}{cccccc}
   Ric_{11}& -w_{34}^2w_{12}&0& -3w_{13}w_{36}w_{34} &0& -2w_{36}^2 w_{13}\\
   -w_{34}^2w_{12} & Ric_{22} &0&0& -w_{13}w_{36}w_{34}&0 \\
    0&0&0&0&0& -w_{36}w_{34}^2 \\
    -3w_{13}w_{36}w_{34}& 0 & 0 & 3w_{36}w_{34}^2 &0& 2w_{36}^2w_{34} \\
    0& -w_{13}w_{36}w_{34}&0&0& w_{36}w_{34}^2 &0 \\
    -2w_{36}^2 w_{13}&0& -w_{36}w_{34}^2 & 2w_{36}^2w_{34} & 0& 2w_{36}^2\\
    \end{array}
\right],
$$
where $Ric_{11} = w_{36}(w_{12}w_{36} + 3w_{13}^2)$, $Ric_{22} = w_{36}(w_{12}w_{36} +w_{13}^2)$,
and scalar curvature $S=\frac{w_{36}^2}{w_{34}^2 w_{12}}.$

Let us formulate the results obtained for the group $G_1$ as follows.

\begin{theorem} \label{G1}
The group $G_1$, which does not admit left-invariant symplectic structures, has multiparameter families of left-invariant semi-K\"{a}hler 2-forms $\w$ with non-degenerate differential $d\w$.
The 3-form $d\w$ invariantly corresponds to the almost para-K\"{a}hler structure $P_{d\w}$, which together with $\w$ forms a semi-para-K\"{a}hler structure $(\w, P_{d\w}, g_{\w})$, where $g_\w(X,Y) = \w(X, P_{d\w}Y)$. The group $G_1$ also has four multiparameter families of semi-K\"{a}hler structures $(\w_k, J, g_k)$, where $J$ is the complex structure (\ref{2}) and $g_k(X,Y) = \w_k(X,JY)$, $k=1,2,3,4$.
\end{theorem}

\subsection{Lie group $G_2$} \label{G2}
This is the second special group that does not admit symplectic structures. The Lie algebra $\g_2$ is defined by the following commutation relations: $[e_1, e_2] = e_4$, $[e_1, e_4] = e_5$, $[e_2, e_4] = e_6$.
This Lie algebra is decomposable: $\g_2 = \h\times \mathbb{R}e_3$, where $\h$ is the subalgebra generated by the vectors $e_1, e_2, e_4, e_5, e_6$.
Let us redesignate the basis vectors so that $\h$ is formed by the first vectors:
$$
e_3 \mapsto e_6, \ e_4 \mapsto e_3, e_5 \mapsto e_4, e_6 \mapsto e_5.
$$
Then the commutation relations will look like this:
$$
[e_1, e_2] = e_3, \ [e_1, e_3] = e_4, \ [e_2, e_3] = e_5.
$$
Lie algebra center: $\mathcal{Z} = \mathbb{R}\{e_4, e_5, e_6\}$.
The Lie algebra $\g_2$ is the semidirect product of the three-dimensional Heisenberg algebra $\h_3$ and the commutative subalgebra $\mathbb{R}^3$:
$$
\g_2 = \h_3 \oplus \mathbb{R}^3 = \mathbb{R}\{e_1, e_3, e_4\}\oplus \mathbb{R}\{e_2, e_5, e_6\}.
$$
%\bowtie
Therefore there is an integrable paracomplex structure on $\g_2$
$$
P = \text{diag}\{+1,-1,+1,+1,-1,-1\}.
$$

\subsubsection{Choice of a non-degenerate semi-K\"{a}hler 2-form.} \label{G2-1}
Let $\w = w_{ij}e^i\wedge e^j$ be an arbitrary 2-form on the Lie algebra $\g_2$.
The exterior derivative $d\w$ is degenerate for any such 2-form $\w$.

The semi-K\"{a}hler property $\w\wedge d\w = 0$ is satisfied under the following conditions:
$$
-w_{14}w_{56} + w_{15}w_{46} -w_{16}w_{45} =0, -w_{24}w_{56} +w_{25}w_{46} -w_{26}w_{45} = 0, -w_{34}w_{56} + w_{35}w_{46} - w_{36}w_{45} = 0.
$$
There is only one solution to this system of equations for a non-degenerate 2-form $\w$.
Then the 2-form $\w$ looks like:
\begin{multline}\label{3}
\w = e^1\wedge  (w_{12} e^2 + w_{13} e^3 + w_{14} e^4 + w_{15} e^5 + w_{16} e^6) + e^2\wedge (w_{23} e^3 + w_{24} e^4 + w_{25} e^5 +\\
+ w_{26} e^6) + e^3\wedge (w_{34} e^4 + w_{35} e^5 + w_{36} e^6).
\end{multline}

Recall that on $\g_2$ there is a natural paracomplex structure $P = \text{diag}\{+1,-1,+1,+1,-1,-1\}$.
The compatible condition $\w(PX,PY) = -\w(X,Y)$ of the form $\w$ and $P$ is satisfied for the following parameter values:
$$
w_{13} = 0, w_{14} = 0, w_{25} = 0, w_{26} = 0, w_{34} = 0.
$$
Then the semi-K\"{a}hler form $\w$ takes the form:
\begin{multline}\label{4}
\w = e^1\wedge  (w_{12} e^2 + w_{15} e^5 + w_{16} e^6) + e^2\wedge (w_{23} e^3 + w_{24}e^4) + e^3\wedge (w_{35} e^5 + w_{36} e^6).
\end{multline}

In this case, we can define the associated metric $g(X,Y) =\w(X,PY)$.
The result is a semi-para-K\"{a}hler structure $(\w, P, g)$ with a Ricci-flat pseudo-Riemannian metric $g$.

\subsubsection{Selecting a compatible complex structure.} \label{G2-2}
It is practically impossible to find a complex structure compatible with the general semi-K\"{a}hler 2-form (\ref{4}).
Therefore, from the family of 2-forms (\ref{4}) we choose the non-degenerate form that is closest to the closed form.
We see that this form vanishes at the center $\mathcal{Z} = \mathbb{R}\{e_4, e_5, e_6\}$ of the Lie algebra.
In our case we have: $de^1 = 0, de^2 = 0, de^6 = 0, de^3 = -e^1\wedge e^2, de^4 = -e^1\wedge e^3, de^5 = -e^2\wedge e^3.$
From these equalities it is easy to see that the terms $w_{12}e^1\wedge e^2$, $w_{13}e^1\wedge e^3$, $w_{23}e^2\wedge e^3$, $w_{16}e^1\wedge e^6$ do not affect the expression of the exterior derivative $d\w$.
Therefore, we will consider $w_{12}$, $w_{13}$, $w_{23}$, and $w_{16}$ to be zero. Then the 2-form $\w$ takes the form:
$$
\w = w_{15} e^1\wedge e^5 + w_{24} e^2\wedge e^4 + e^3\wedge (w_{35} e^5 + w_{36} e^6).
$$
This form is non-degenerate for $w_{15}w_{24}w_{36} \ne 0$.
This condition does not contain $w_{35}$.
Therefore, we will also assume that $w_{35}=0$. We obtain the following non-degenerate semi-K\"{a}hler 2-form
$$
\w = w_{15} e^1\wedge e^5 + w_{24} e^2\wedge e^4 + w_{36} e^3\wedge e^6.
$$
By renormalizing the basis vectors of the Lie algebra, this form can be reduced to the form:
\begin{equation} \label{5}
\w_0 = e^1\wedge e^5 -e^2\wedge e^4 -e^3\wedge e^6.
\end{equation}

Now we will look for a nilpotent complex structure $J = (\psi_{ij})$ compatible with the resulting semi-K\"{a}hler form (\ref{5}).
To do this, you need to find $J$ that satisfies three conditions: $\w_0(X,JY) + \w_0(JX,Y) = 0$, $J^2 = -Id$ and $N_J = 0$, which are expressed by the system of equations (\ref{0}).

The compatible condition $\w_0(X,JY) + \w_0(JX,Y) = 0$ is satisfied for the following parameter values:
$$
\psi_{52} = -\psi_{41}, \psi_{61} = -\psi_{53}, \psi_{54} = \psi_{21}, \psi_{55} = -\psi_{11}, \psi_{56} = \psi_{31}, \psi_{62} = \psi_{43}, \psi_{44} = -\psi_{22}, \psi_{45} = \psi_{12},
$$
$$
\psi_{46} = -\psi_{32}, \psi_{64} = -\psi_{23}, \psi_{65} = \psi_{13}, \psi_{66} = -\psi_{33}, \psi_{25} = -\psi_{14}, \psi_{26} = \psi_{34}, \psi_{16} = -\psi_{35}.
$$
Consider an increasing sequence of $J$-invariant ideals:
$$
\mathfrak{a}_0(J) = 0, \mathfrak{a}_s(J) = \{X \in \g | [X, \g] \subset \mathfrak{a}_{s-1}(J) \mbox{ and } [JX, \g] \subset \mathfrak{a}_{s-1} (J)\}, \ s \ge 1.
$$
In particular, the ideal $\mathfrak{a}_1(J)$ lies in the center and has dimension at least two.
In our case, we have $\mathfrak{a}_1(J)\subset \mathcal{Z} = \mathbb{R}\{e_4, e_5, e_6\}$ and has dimension 2.
From the commutation relations $[e_1, e_2] = e_3, \, [e_1, e_3] = e_4, \, [e_2, e_3] = e_5$ it follows that $\mathfrak{a}_2(J) = \mathbb{R}\{e_3, e_4, e_5, e_6\}$ and $\mathfrak{a}_3(J) = \g_2$.
Taking into account the compatible property, we obtain the following form of the matrix of the nilpotent almost complex structure $J$:
$$
J =\left(
  \begin{array}{cccccc}
    \psi_{11} & \psi_{12} &0&0&0&0 \\
    \psi_{21} & \psi_{22} &0&0&0&0 \\
    \psi_{31} & \psi_{32} & \psi_{33} &0&0& \psi_{36} \\
    \psi_{41} & \psi_{42} & \psi_{43} & -\psi_{22} & \psi_{12} & -\psi_{32} \\
    \psi_{51} & \psi_{41} & \psi_{53} & \psi_{21} & -\psi_{11} & \psi_{31}\\
    -\psi_{53} & \psi_{43} & \psi_{63} &0&0& -\psi_{33} \\
  \end{array}
\right),
$$
with the following parameter values:
$$
\psi_{13} = 0, \psi_{14} = 0, \psi_{15} = 0, \psi_{35} = 0, \psi_{23} = 0, \psi_{24} = 0, \psi_{34} = 0.
$$
Now we solve two systems of equations: $J^2 = -Id$ and $N_J = 0$.
As a result, we obtain the following left-invariant complex structure:
\begin{equation} \label{6}
J =\left(
  \begin{array}{cccccc}
    \psi_{11} & \psi_{12} &0&0&0&0 \\
    -\frac{w_{11}^2+1}{\psi_{12}} & -\psi_{11} &0&0&0&0 \\
    J^3_1 & \psi_{32} & \psi_{33} &0&0& -\frac{w_{33}^2+1}{\psi_{63}} \\
    \psi_{41} & \psi_{42} & \psi_{43} & \psi_{11} & \psi_{12} & -\psi_{32} \\
    J^5_1 & -\psi_{41} & J^5_3  & -\frac{w_{11}^2+1}{\psi_{12}} & -\psi_{11} & J^5_6\\
    J^6_1  & \psi_{43} & \psi_{63} &0&0& -\psi_{33} \\
  \end{array}
\right),
\end{equation}
where $J_n^k$ are rational functions of parameters $\psi_{ij}$:
$$
J_1^3=J_6^5 =\frac{\psi_{11} \psi_{32} \psi_{63} -\psi_{32}\psi_{33} \psi_{63}+\psi_{33}^2 \psi_{43}+\psi_{43}}{\psi_{12} \psi_{63}},
$$
$$
J_1^5 =\frac{\psi_{11}^2 \psi_{42}\psi_{63} -2\psi_{11}\psi_{12}\psi_{41} \psi_{63} -\psi_{32}^2 \psi_{63}^2 +2\psi_{32} \psi_{33} \psi_{43} \psi_{63} -\psi_{33}^2 \psi_{43}^2 +\psi_{42}\psi_{63} -\psi_{43}^2}{\psi_{12}^2 \psi_{63} },
$$
$$
J_1^6=-J_3^5 =\frac{\psi_{11} \psi_{43} -\psi_{32} \psi_{63} +\psi_{33} \psi_{43}}{\psi_{12} }.
$$

Let $g_J(X,Y) = \w(X,JY)$ be the associated metric. We obtain a family of semi-K\"{a}hler structures $(\w, J, g_J)$ with a Ricci-flat pseudo-Riemannian metric $g_J$.

\begin{theorem} \label{G2}
The group $G_2$, which does not admit left-invariant symplectic structures, admits a multiparameter family of left-invariant semi-para-K\"{a}hler structures $(\w,P,g)$ with the semi-K\"{a}hler 2-form (\ref{4}), the integrable paracomplex structure $P = diag\{+1,-1,+1,+1,-1, -1\}$ and the pseudo-Riemannian Ricci-flat metric $g(X,Y) =\w(X,PY)$. The group $G_2$ also admits a multiparameter family of semi-K\"{a}hler structures $(\w_0,J,g_J)$ with fundamental form (\ref{5}), nilpotent complex structure (\ref{6}) and associated pseudo-Riemannian Ricci-flat metric $g_J(X,Y) =\w_0(X,JY)$.
\end{theorem}

\subsection{Lie group $G_3$} \label{G3}
%This is the third special group, which does not admit symplectic structures, but has complex ones.
The Lie algebra $\g_3$ is defined by the following commutation relations:
$[e_1, e_2] = e_6$, $[e_3, e_4] = e_6$. Let us rename $e_6$ to $e_5$, then
$$
[e_1, e_2] = e_5, [e_3, e_4] = e_5.
$$
The Lie algebra is decomposable: $\g_3 = \h\times \mathbb{R}e_6$. Here the algebra $\h$ is the central extension of $\mathbb{R}^4$ with the standard symplectic form $\w = e^1\wedge e^2 + e^3\wedge e^4$ and $\mathbb{R}e_5$.

Recall that if there is a symplectic Lie algebra $(\n,\w)$, then the central extension $\h = \n\times_\w \mathbb{R}$ is a Lie algebra in which the Lie brackets are defined as follows:
$$
[X,\xi]_\h = 0, \ [X,Y]_\h = [X,Y]_\n + \w(X,Y)\xi,
$$
for any $X,Y\in \n$, where $\xi = d/dt$ is a unit vector from $\mathbb{R}$.

In our case, $\h = \mathbb{R}^4\times_\w \mathbb{R}e_5$, $\w = e^1\wedge e^2 + e^3\wedge e^4$ and $\xi = e_5$.
Since the form $\w$ is non-degenerate, $\h$ is a contact Lie algebra \cite{Diatta} with contact form $\eta = e^5$ and Reeb field $\xi =e_5$.

Let us define the affinor $\varphi$ on $\h$ in the following natural way:
$$
\varphi(e_1) = e_2, \quad  \varphi(e_2) =-e_1, \quad \varphi(e_3) =e_4, \quad \varphi(e_4) =-e_3,\quad \varphi(e_5) = 0.
$$
The associated metric for the contact structure $\eta$ is completely determined by the affinor $\varphi$ by the formula:
$$
g(X,Y) = d\eta(\varphi X,Y) + \eta(X)\eta(Y).
$$
Recall that a contact metric structure $(\eta,\xi,\varphi,g)$ on a manifold $M$ is called a Sasaki structure if the almost complex structure $J$ on $M\times \mathbb{R}$ defined by the formula is integrable
$$
J(X,f\partial) = (\varphi X -f\xi, \eta(X)\partial).
$$
Here the tangent vector to $M\times \mathbb{R}$ is represented as a pair $(X,f\partial)$, where $X$ is the tangent vector to $M$ and $f\partial$ is the tangent vector to $\mathbb{R}$, $\partial$ is the basis vector on $\mathbb{R}$.
In our case, $M = G_3$, $\xi = e_5$, $M\times \mathbb{R} = G_3 \times \mathbb{R}e_6$. Then we define an almost complex structure on $\g_3 = \h\times \mathbb{R}e_6$ by the formula:
$$
J(X,fe_6) = (\varphi X -fe_5, \eta(X)e_6).
$$
In this case
$$
J(e_i) = J(e_i,0\cdot e_6) = (\varphi e_i, \eta(e_i)e_6) = (\varphi e_i, 0) = \varphi (e_i),\quad i = 1,2,3,4,
$$
$$
J(e_5) = J(e_5, 0\cdot e_6) = (\varphi e_5, \eta(e_5)e_6) = (0, e_6) = e_6,
$$
$$
J(e_6) = J(0, e_6) = (0 - e_5, \eta(0)e_6) = (-e_5, 0) = -e_5.
$$
Calculations show that this almost complex structure $J$ is integrable. Therefore $(\eta,\xi,\varphi,g)$ is a Sasaki structure on $\h$.

Let us define a 2-form on $\g_3 = \h\times \mathbb{R}e_6$
\begin{equation} \label{7}
\w = e^1\wedge e^2 + e^3\wedge e^4 + e^5\wedge e^6.
\end{equation}
and the associated metric $g(X,Y) =\w(X,JY)$. We obtain a Hermitian structure $(\w, J, g)$ of scalar curvature $S =-1$ and with a Ricci tensor of the form:
$$
Ric =\frac 12 \left(
  \begin{array}{cccccc}
    -1&0&0&0&0&0 \\
    0&-1&0&0&0&0 \\
    0&0&-1&0&0&0 \\
    0&0&1&-1&0&0 \\
    0&0&0&0&2&0 \\
    0&0&0&0&0&0 \\
  \end{array}
\right).
$$

\subsubsection{Semi-K\"{a}hler structures.}
Unfortunately, the 2-form (\ref{7}) does not have the semi-K\"{a}hler property: $\w\wedge d\w \ne 0$.
However, the following 2-form
\begin{equation} \label{8}
\w = e^1\wedge e^2 - e^3\wedge e^4 + e^5\wedge e^6.
\end{equation}
is semi-K\"{a}hler.
Let us find a nilpotent complex structure compatible with the form (\ref{8}).
We take the general almost complex structure $J = (\psi_{ij})$ and require the invariance of the center $\mathcal{Z}$, the compatible condition $\w(X,JY) + \w(JX,Y) = 0$,  the property $J^2 = -Id$, and the integrability condition $N_J = 0$.

The invariance of the center $\mathcal{Z} = \mathbb{R}\{e_5,e_6\}$ occurs for the following parameter values:
$$
\psi_{15} = 0, \psi_{16} = 0, \psi_{25} = 0, \psi_{26} = 0, \psi_{35} = 0, \psi_{36} = 0, \psi_{45} = 0, \psi_{46} = 0.
$$
Compatible condition is satisfied under the following conditions on the parameters:
$$
\psi_{22} = -\psi_{11}, \psi_{41} = -\psi_{23}, \psi_{24} = \psi_{31}, \psi_{61} = 0, \psi_{51} = 0, \psi_{42} = \psi_{13}, \psi_{14} = -\psi_{32}, \psi_{62} = 0, \psi_{52} = 0,
$$
$$
\psi_{44} = -\psi_{33}, \psi_{66} = -\psi_{55}, \psi_{63} = 0, \psi_{64} = 0, \psi_{53} = 0, \psi_{54} = 0.
$$
As a result, we obtain the following form of the matrix of the almost complex structure $J$:
$$
J =\left(
  \begin{array}{cccccc}
-\psi_{22}&\psi_{12}&\psi_{13}&\psi_{14} &0&0 \\
\psi_{21} & \psi_{22} & \psi_{23} & \psi_{24} &0&0 \\
-\psi_{24} & \psi_{14} & -\psi_{33} & \psi_{34} &0&0 \\
\psi_{23} & -\psi_{13} & \psi_{43} & \psi_{33} &0&0 \\
0&0&0&0& -\psi_{66} & \psi_{56} \\
0&0&0&0& \psi_{65} & \psi_{66} \\
  \end{array}
\right).
$$

We require that the condition $J^2 = -Id$ be satisfied and that the Nijenhuis tensor be equal to zero. We get:
\begin{equation} \label{9}
J =\left(
  \begin{array}{cccccc}
\psi_{11} & \psi_{12} &0&0&0&0 \\
-\frac{\psi_{11}^2+1}{\psi_{12}} & -\psi_{11} &0&0&0&0 \\
0&0& \psi_{33} & \psi_{34} &0&0 \\
0&0& -\frac{\psi_{33}^2+1}{\psi_{34}} & -\psi_{33} &0&0 \\
0&0&0&0& \psi_{55} & \psi_{56} \\
0&0&0&0& -\frac{\psi_{55}^2+1}{\psi_{56}} & -\psi_{55} \\
  \end{array}
\right).
\end{equation}
Let $g_J(X,Y) =\w(X,JY)$ be the associated metric. We obtain a family of semi-K\"{a}hler structures $(\w,J,g_J)$ with the Ricci trnsor
$$
Ric := \frac{\psi_{55}^2+1}{2\psi_{56}}
\left( \begin {array}{cccccc}
-{\frac {{\psi_{11}}^{2}+1}{\psi_{12}}} & -\psi_{11} & 0 & 0 & 0 & 0\\
-\psi_{11} & -\psi_{12} & 0 & 0 & 0 & 0 \\
0 & 0 & {\frac {{\psi_{33}}^{2}+1}{\psi_{34}}} & \psi_{33} & 0 & 0\\
0 & 0 & \psi_{33} & \psi_{34} & 0 & 0 \\
0 & 0 & 0 & 0 & 2\,{\frac {{\psi_{55}}^{2}+1}{\psi_{56}}} & 2\,\psi_{55}\\
0 & 0 & 0 & 0 & 2\, \psi_{55} & 2\,{\frac {\psi_{56}\,{\psi_{55}}^{2}}{{\psi_{55}}^{2}+1}}
\end {array} \right)
$$
and scalar curvature $S=\frac{\psi_{55}^2+1}{2\psi_{56}}$.

\subsubsection{Semi-para-K\"{a}hler structures.}
The Lie algebra $\g_3$ is the semidirect product of two commutative subalgebras generated by the vectors $e_1, e_3, e_5$ and $e_2, e_4, e_6$, respectively. Therefore there is an integrable paracomplex structure on $\g_3$:
$$
P = \text{diag}\{+1,-1,+1, -1,+1,-1\}.
$$
The compatibley condition $\w(PX,PY) = -\w(X,Y)$ with the 2-form (\ref{8}) is obviously satisfied. Let $g(X,Y) =\w(X,PY)$ be the associated metric. We obtain the semi-para-K\"{a}hler structure $(\w,P,g)$ with a Ricci-flat pseudo-Riemannian metric $g$.

Let us formulate the results obtained for the group $G_3$ in the following form.

\begin{theorem} \label{G3}
The group $G_3$ that does not admit left-invariant symplectic structures has a multiparameter family (\ref{9}) of left-invariant nilpotent complex structures $J$ compatible with the semi-K\"{a}hler 2-form (\ref{8}). In this case, the triple $(\w,J,g)$, where $g(X,Y) =\w(X,JY)$, defines a semi-K\"{a}hler structure of nonzero Ricci curvature. The group $G_3$ also admits a left-invariant semi-para-K\"{a}hler Ricci-flat structure.
\end{theorem}

\end{document}